\documentclass{elsarticle}
\usepackage{mathptmx}       
\usepackage{helvet}         
\usepackage{courier}        
\usepackage{type1cm}        
\usepackage{graphicx}        

\usepackage{amsmath}			
\usepackage{amsfonts}			
\usepackage{amsthm}			
\newtheorem{theorem}{Theorem}

\newtheorem{mycor}{Corollary}
\usepackage{silence}
\usepackage{pgfplots}
\pgfplotsset{compat=newest,every axis/.append style={legend style={font=\tiny}}}
\usepackage{pgfplotstable}

\usepackage{hyperref}
\usepackage[capitalize]{cleveref}
\crefname{assumption}{assumption}{assumptions}
\creflabelformat{assumption}{(#2#1#3)}
\crefname{inequality}{inequality}{inequalities}
\creflabelformat{inequality}{(#2#1#3)}
\crefname{equation}{}{}

\usepackage{todonotes}

\newcommand{\levy}{L\'{e}vy}

\renewcommand{\d}{\mbox{d}}
\newcommand{\W}{\mbox{W}}

\newcommand{\E}{\mathbb{E}}
\newcommand{\Prob}{\mathbb{P}}

\newcommand{\bigO}{\mathcal{O}}
\begin{document}
\title{Weak Antithetic MLMC Estimation of SDEs with the Milstein scheme for Low-Dimensional Wiener Processes}

\author[sdu]{Kristian Debrabant}
\ead{debrabant@imada.sdu.dk}
\author[sdu,iut]{Azadeh Ghasemifard}
\ead{azadeh.ghasemi@math.iut.ac.ir}
\author[sdu]{Nicky C.\ Mattsson\texorpdfstring{\corref{cor1}}{}}
\ead{mattsson@imada.sdu.dk}

\cortext[cor1]{Corresponding author}
\address[sdu]{Department of Mathematics and Computer Science, University of Southern Denmark, Odense, Denmark}
\address[iut]{Department of Mathematical Sciences, Isfahan University of Technology, P.O. Box 8415683111, Isfahan, Iran}

\begin{keyword}
multilevel Monte-Carlo \sep stochastic differential equation \sep weak approximation schemes \sep Milstein scheme
\MSC[2010] 65C30 \sep 65C05
\\
~\\
\copyright 2018. This manuscript version is made available under the CC-BY-NC-ND 4.0 license \url{http://creativecommons.org/licenses/by-nc-nd/4.0/}
\end{keyword}

\begin{abstract}
In this paper, we implement a weak Milstein Scheme to simulate low-dimensional stochastic differential equations (SDEs). We prove that combining the antithetic multilevel Monte-Carlo (MLMC) estimator introduced by Giles and Szpruch with the MLMC approach for weak SDE approximation methods by Belomestny and Nagapet\-yan, we can achieve a quadratic computational complexity in the inverse of the Root Mean Square Error (RMSE) when estimating expected values of smooth functionals of SDE solutions, without simulating \levy\ areas and without requiring any strong convergence of the underlying SDE approximation method. By using appropriate discrete variables this approach allows us to calculate the expectation on the coarsest level of resolution by enumeration, which, for low-dimensional problems, results in a reduced computational effort compared to standard MLMC sampling. These theoretical results are also confirmed by a numerical experiment.
\end{abstract}

\maketitle

\section{Introduction}
In this paper we inquire to efficiently estimate $\E[f(X(T))]$, where $f:\mathbb{R}^{d}\rightarrow\mathbb{R}$ is a sufficiently smooth Lipschitz continuous function and $X(t)$ follows the stochastic differential equation (SDE)
\begin{equation} \label{equ:diffusion}
X(t) = X_0 + \int_{t_0}^{t}\mu(X(s)) \d s + \int_{t_0}^{t} \sigma(X(s))\d \W (s),
\end{equation}
with $X_0 \in \mathbb{R}^d$ being a known initial condition, $\mu \in C^2(\mathbb{R}^d,\mathbb{R}^d)$, $\sigma \in C^2 (\mathbb{R}^{d \times m},\mathbb{R}^d)$ and $\W(t)$ being an $m$-dimensional Wiener process.

The above problem typically occurs in option pricing, where $f$ is the payoff one is given of the underlying asset $X(t)$ \cite{black1973pricing,kloeden1992numerical,giles2008multilevel,gilesXXaom}.

A decade ago, Giles \cite{giles2008multilevel} introduced the Multilevel Monte-Carlo (MLMC) method to estimate $\E[f(X(t))]$, which combines approximations with varying number of paths on different levels of resolution and results often in a significant reduction of the computational complexity over standard Monte-Carlo. Proving the requirements to obtain this reduction was initially usually based on the strong convergence properties of the underlying SDE approximation method. Later \cite{giles2014antithetic} and \cite{belomestny2014multilevel} managed to circumvent this strong convergence requirement by introducing antithetic paths and a coupling between the levels respectively.

In this paper we combine these efforts to construct a method which has a computational complexity of optimal order, without requiring any strong convergence of the underlying SDE approximation method. This allows further to use enumeration to calculate the expectation on the coarsest level of resolution, resulting in a reduced computational complexity in case of low dimensional Wiener processes.

We continue this paper by shortly describing the standard MLMC and the two extensions in \cref{sec:MLMC}, showing that a combination of these two efforts is possible. Hereafter we introduce in \cref{sec:enumeration} the idea of enumeration, which has the potential to significantly reduce the number of paths on the initial level. We then show the practical usefulness of this approach when pricing a basket option in \cref{sec:results}.

\section{Multilevel Monte-Carlo Estimation}\label{sec:MLMC}
Giles introduced in \cite{giles2008multilevel} the MLMC approach which greatly improves on the standard Monte-Carlo sampling. Let $P=f(X(T))$ and $\Delta_l = 2^{-l}T$ be the stepsize used to obtain the approximation $\hat{P}_l$ of $f(X(T))$ at level $l$, $l=0,\ldots, L$. Then instead of applying standard Monte-Carlo to $\E(\hat{P}_L)$ we apply Monte-Carlo to the telescopic sum $\E(\hat{P}_L) = \E(\hat{P}_0^{f}) + \sum_{l = 1}^{L} \E(\hat{P}^{f}_l  - \hat{P}_{l-1}^{c})$ with different numbers $M_l$ of paths and $\hat{P}^{f}_l$ and $\hat{P}^{c}_l$ such that $\E(\hat{P}^{f}_l)=\E(\hat{P}^{c}_l)=\E(\hat{P}_l)$, resulting in the estimate
\begin{equation*}
\hat{Y} = \sum_{l=0}^{L} \hat{Y}_l, \quad \mbox{where } \hat{Y}_l =\begin{cases}
M_l^{-1} \sum_{i=1}^{M_l} (\hat{P}_l^{f,i} - \hat{P}_{l-1}^{c,i}), & \mbox{if } l\neq 0 \\
M_0^{-1} \sum_{i=1}^{M_0} \hat{P}_0^{f,i}, & \mbox{if } l= 0
\end{cases}
\end{equation*}
where the upper index $i$ denotes the $i$-th realization.

Giles \cite{giles2008multilevel} proved that under some conditions the computational complexity, $C$, of this scheme is directly related to the order $\beta$ of variance reduction, $\mbox{var}(\hat{Y}_l)=\bigO(\Delta_l^{\beta})$, and the Mean Square Error (MSE), $\epsilon^2= \E(\hat{Y}-\E(P))^2$ through
\begin{equation} \label{equ:compcost}
C = \begin{cases}
\bigO(\epsilon^{-2}), & \mbox{if } \beta >1 \\
\bigO(\epsilon^{-2} (\log\epsilon)^2), & \mbox{if } \beta =1.
\end{cases}
\end{equation}
As $\beta$ is bounded from below by twice the order of strong convergence, initially SDE approximation methods with a strong order at least 0.5, ideally $\geq1$, seemed necessary. Recently there has been some research trying to circumvent this strong order requirement. We will in the following describe and combine these efforts to obtain a method which solves \cref{equ:diffusion} with $\beta=2$ without being strong convergent at all.

\subsection{Weak MLMC}
In \cite{belomestny2014multilevel} Belomestny and Nagapetyan introduced the Weak MLMC method, proving that MLMC based on the weak Euler scheme maintains $C=\bigO(\epsilon^{-2} (\log\epsilon)^2)$ when doing the right coupling of levels. Concretely, for an arbitrary level $l$, let $\xi^f_{l,i}$ and $\xi^c_{l,i}$, $i=1,\dots,2^l$, be, possibly approximate, Wiener increments with variance $\Delta_{l}$, used in the approximation of $\hat{P}^{f}_l$ and $\hat{P}^{c}_l$, respectively. Then we require
\begin{equation} \label{equ:coupling}
\mathcal{R} = \xi^c_{l-1,i} - \xi^f_{l,2i-1} - \xi^f_{l,2i}, \quad \xi^c_{l-1,i} \stackrel{D}= \xi^f_{l-1,i}
\end{equation}
with sufficiently small $\mathcal{R}$. We will here only consider $\mathcal{R}=0$.

\Cref{equ:coupling} couples the random variables on a single level together and to the lower levels. One picks an initial distribution for the random variables $\xi^f_{L,i}$, which then by \cref{equ:coupling} generates the distributions of the random variables on the lower levels. For example, the choice $\xi^f_{L,i} \sim \mbox{N}(0,\sqrt{\Delta_{L}})^m$ yields $\xi^f_{l,i} \sim \mbox{N}(0,\sqrt{\Delta_{l}})^m$, the classical MLMC, while for the two-point distribution approximation with $\Prob(\xi^f_{L,i,j} = \pm \sqrt{\Delta_{L}})=\frac{1}{2}$, $j=1,\dots,m$, we obtain
\begin{equation} \label{equ:weakBrownian}
\begin{gathered}
 \xi_{l,i}^f \sim \left( \mbox{Bin}\left( 2^{L-l},0.5\right) - 2^{L-l-1} \right) \cdot 2\sqrt{\Delta_L}.
\end{gathered}
\end{equation}
It is in \cite[Section 4.1]{belomestny2014multilevel} discussed how to implement the generation of binomial random numbers efficiently.

\subsection{Weak Antithetic MLMC}
One large problem with the strong convergence requirement is that a first order strong convergent method is not readily available for multiple dimensions, due to troublesome terms known as \levy -areas, that are in general computationally quite demanding.

Consider (for $\alpha\in\{c,f\}$) the multidimensional Milstein scheme for SDEs
\begin{equation} \label{equ:Milstein}
\begin{aligned}
X^{\alpha}_{l,i}=X^{\alpha}_{l,i-1}&+\mu(X_{l,i-1}^{\alpha})\Delta_{l} +\sum_{j=1}^{m}\sigma_{j}(X^{\alpha}_{l,i-1}) \xi^{\alpha}_{l,i-1}\\
&+\sum_{j,k=1}^{m}\sigma_j^{'}(X^{\alpha}_{l,i-1})\sigma_k(X^{\alpha}_{l,i-1})(\xi^{\alpha}_{l,i-1,j} \xi^{\alpha}_{l,i-1,k}-\Omega_{jk}\Delta_{l}-A^{\alpha}_{jk,i-1}),
\end{aligned}
\end{equation}
with $\Omega$ being the correlation matrix and $A^{\alpha}_{jl,i-1}$ being the \levy -areas as defined in \cite{giles2014antithetic}, and $\xi^{\alpha}_{l,i-1,k}$ being the $k$-th component of the Wiener increment with variance $\Delta_{l}$.

Giles and Szpruch \cite{giles2014antithetic} showed that if one also defines an antithetic path $X^a_{l,i}$, obtained by switching $\xi^f_{l,2i-1,j}$ with $\xi^f_{l,2i,j}$ for all $i$, and chooses $\hat{P}^{c}_l=f(X^{c}_{l,T/\Delta_l})$ and $\hat{P}^{f}_l= 1/2(f(X^f_{l,T/\Delta_l}) + f(X^a_{l,T/\Delta_l}))$, one can ignore the \levy -area while maintaining $\beta=2$ for smooth payoff functions $f$.

We now extend this result to the scheme where $\xi^c_{l,i-1},\xi^f_{l,i-1}$ are approximate Wiener increments as given in \cref{equ:coupling,equ:weakBrownian}.

\begin{theorem}\label{theorem:MilsteinVarRed}
Using the Milstein method without \levy\ areas,
\begin{equation} \label{equ:WeakMilstein}
\begin{aligned}
X^{\alpha}_{l,i}=X^{\alpha}_{l,i-1}&+\mu(X^{\alpha}_{l,i-1})\Delta_{l} +\sum_{j=1}^{m}\sigma_{j}(X^{\alpha}_{l,i-1}) \xi^{\alpha}_{l,i-1}\\
&+\sum_{j,k=1}^{m}\sigma_j^{'}(X^{\alpha}_{l,i-1})\sigma_k(X^{\alpha}_{l,i-1})(\xi^{\alpha}_{l,i-1,j} \xi^{\alpha}_{l,i-1,k}-\Omega_{jk}\Delta_{l}),
\end{aligned}
\end{equation}
with $\alpha\in\{c,f\}$ and $\xi^{\alpha}_{l,i}$ being the weak approximation given by \cref{equ:coupling,equ:weakBrownian}, together with the antithetic MLMC to estimate $\E(f(X(T))$, yields second order variance reduction for Lipschitz continuous payoff functions $f$ in $C^2(\mathbb{R}^d,\mathbb{R})$.
\begin{proof}
This proof is similar to the one given in \cite{giles2014antithetic}, with the change that as \cite[Lemma 4.6]{giles2014antithetic} requires strong convergence to bound $\Vert X^f_{l,i} - X^a_{l,i} \Vert$, to circumvent this we use an adapted version of \cite[Theorem 4.10]{giles2014antithetic}.
\end{proof}
\end{theorem}

\begin{mycor}[Computational Complexity]
Under the conditions of \cref{theorem:MilsteinVarRed}, using the Milstein scheme without Levy areas together with binomial random numbers and the antithetic MLMC results in a computational complexity of $C = \bigO(\epsilon^{-2})$.
\begin{proof}
Combining \cref{theorem:MilsteinVarRed} with \cite[Theorem 3.1]{giles2008multilevel} yields the desired result.
\end{proof}
\end{mycor}

Note: For the European option, which is not in $C^2$, we have the more conservative bound $\beta>1.5$ \cite[Theorem 5.2]{giles2014antithetic}. However, this still gives $C = \bigO(\epsilon^{-2})$ by \cref{equ:compcost}.

\section{The Enumeration Idea}\label{sec:enumeration}
Using discrete variables, e.g.\ \cref{equ:weakBrownian}, instead of the normally distributed random variables representing Wiener increments allow us to calculate the expectation exactly by calculating all possible outcomes and scaling by their probability. We will in this paper solely consider \cref{equ:weakBrownian}.

Let $\bar{x} = (x_1, \ldots , x_m)^T$ be a sample of the random variable $\xi^f_{l,i} = (\xi^f_{l,i,1}, \ldots ,\xi^f_{l,i,m})$. As the different approximate Wiener increments are independent the probability of each collection of outcomes is
\begin{equation*}
\Prob(\xi^f_{l,i} = \bar{x}) = \prod_{k=1}^{m} \Prob(\xi^{f}_{l,i,k}=x_k), \quad \mbox{with }\Prob(\xi^{f}_{l,i,k}=x_k) = \begin{pmatrix} 2^{L-l} \\ \frac{x_k}{2\sqrt{\Delta_L}} + 2^{L-l-1} \end{pmatrix} 0.5^{2^{L-l}}.
\end{equation*}
As each "Binomial" distribution has $2^{L-l}+1$ different outcomes, there are $2^{l}$ independent steps and $m$ independent dimensions, the number of possible paths on level $l$ is  given by $(2^{L-l}+1)^{2^l m}$. Due to this, we apply enumeration only on level $l=0$, where the number of paths is largest in the MLMC. The number of possible paths is then given by $\hat{M}=(2^{L}+1)^{m}$.

In \cref{tbl:paths} we see that enumeration is almost unconditionally useful in the cases of $m=1$ and $m=2$, whereas for $m=3$ and $m=4$ it is still useful if the number of levels can be kept low, which often is the case, see e.\,g.\ the example in \cref{sec:results}. For $m>4$ enumeration is only rarely beneficial, unless the number of levels needed is low and the variance extraordinarily large.
\begin{table}[ht!]
\begin{tabular}{|c|c|c|c|c|}
\hline
$\hat{M}$ & $m=1$ & $m=2$ & $m=3$ & $m=4$ \\
\hline
$L=8$ & $257$ & $66,049$ & $16,974,593$ & $4,362,470,401$ \\
\hline
$L=9$ & $513$ & $263,169$ & $135,005,697$ & $-$ \\
\hline
$L=10$ & $1,025$ & $1,050,625$ & $1,076,890,625$ & $-$ \\
\hline
$L=11$ & $2,049$ & $4,198,401$ & $8,602,523,649$ & $-$ \\
\hline
\end{tabular}
\caption{Number of simulations $\hat{M}$ needed to calculate $\E(\hat{P}_0)$ as a function of the number of dimensions for the Wiener process and the maximum level. Dashes mark places where there are more than $1e10$ possibilities.} \label{tbl:paths}
\end{table}

\section{Numerical Example} \label{sec:results}
In this section we price the basket option
\begin{equation}
f(X_1(T),X_2(T)) = e^{-rT} \max\left( X_1(T)+X_2(T) - K,0 \right),
\end{equation}
using the following linear SDE to simulate the underlying asset,
\begin{align*}
\d X_1(t) &= r X_1(t) \d t + \sigma_1 X_1(t) \d \W_1(t) + \sigma_2 (X_1(t)+X_2(t))\d \W_2(t), \\
\d X_2(t) &= r X_2(t) \d t + \sigma_3 X_2(t) \d \W_1(t) + \sigma_4 X_2(t) \d \W_2(t).
\end{align*}
We use the parameters $T=1$, $K=2$, $X_k(0)=1$, $r = 0.2$, $\sigma_1 = 0.05$, $\sigma_2 = 0.1$, $\sigma_3 = 0.15$, $\sigma_4 =0.2$ and compare three different MLMC methods:
\begin{itemize}
\item the Euler MLMC method with binomial random variables introduced by Belomestny and Nagapetyan \cite{belomestny2014multilevel} (Euler),
\item the antithetic Milstein MLMC with normal random variables introduced by Giles and Szpruch \cite{giles2014antithetic} (Normal),
\item the antithetic Milstein MLMC with binomial random variables and enumeration proposed in this article (Binomial).
\end{itemize}
We will use the name in parenthesis as classifier in the legend. For all of them we approximate the variance reduction and convergence for root mean square error $RMSE=1e-5$, as well as the number of paths on $l=0$ and the computational complexity for $RMSE \in \{5e-4,2e-4,1e-4,5e-5,2e-5,1e-5\}$. We have repeated this $100$ times and report in \cref{fig:EuropeanCall} the average of these results together with their $95\%$ confidence intervals, often vanishingly small.

For the Binomial and Euler method we approximate the number of levels from a simulation with $RMSE=5e-4$ and the assumption of linear convergence and allow the code to increase the levels if necessary, this is included in the computational complexity. For the Normal method we can reuse already calculated paths when $L$ increases, thus here we just start with $L=1$.

\pgfplotstableread[col sep=comma]{Figures/StrongMilstein/variance1.csv}\StrongVar
\pgfplotstableread[col sep=comma]{Figures/StrongMilstein/convergence1.csv}\StrongCon
\pgfplotstableread[col sep=comma]{Figures/StrongMilstein/compcost1.csv}\StrongComp
\pgfplotstableread[col sep=comma]{Figures/StrongMilstein/paths1.csv}\StrongPath

\pgfplotstableread[col sep=comma]{Figures/WeakMilstein/variance1.csv}\WeakVar
\pgfplotstableread[col sep=comma]{Figures/WeakMilstein/convergence1.csv}\WeakCon
\pgfplotstableread[col sep=comma]{Figures/WeakMilstein/compcost1.csv}\WeakComp
\pgfplotstableread[col sep=comma]{Figures/WeakMilstein/paths1.csv}\WeakPath

\pgfplotstableread[col sep=comma]{Figures/Belomestny/variance1.csv}\BeloVar
\pgfplotstableread[col sep=comma]{Figures/Belomestny/convergence1.csv}\BeloCon
\pgfplotstableread[col sep=comma]{Figures/Belomestny/compcost1.csv}\BeloComp
\pgfplotstableread[col sep=comma]{Figures/Belomestny/paths1.csv}\BeloPath

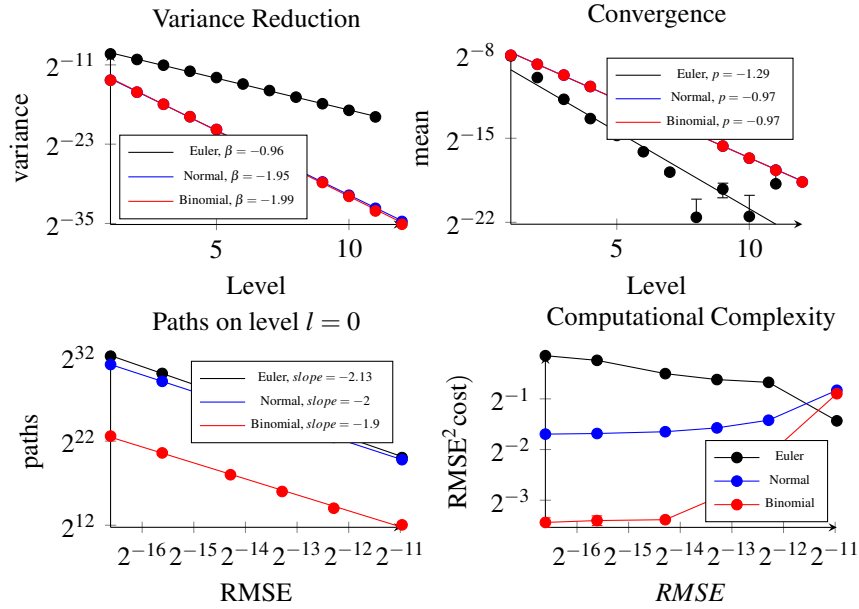
\begin{figure}[ht!]
\begin{tikzpicture}
\begin{axis}[%
    name=plot1,
    width=0.45\textwidth,
    height=0.2\textheight,
    ymode = log,
    log basis y={2},
    axis x line=bottom,
    axis y line=left,
    title=Variance Reduction,
    xlabel=Level,
    ylabel={variance},
    legend pos=south west]

	\addplot[only marks, black, forget plot,error bars/.cd, y dir=both, y explicit] table[x={level}, y={variance}, y error={err}]{\BeloVar};
    \addplot[mark=none, black] table[x={level}, y={create col/linear regression={y=variance}}]{\BeloVar};
    \xdef\slopeBelomestny{\pgfplotstableregressiona}
	\addlegendentry{Euler, $\beta=\pgfmathprintnumber{\slopeBelomestny}$}

    \addplot[only marks, blue, forget plot,error bars/.cd, y dir=both, y explicit] table[x={level}, y={variance}, y error={err}]{\StrongVar};
    \addplot[mark=none, blue] table[x={level}, y={create col/linear regression={y=variance}}]{\StrongVar};
    \xdef\slopeStrong{\pgfplotstableregressiona}
	\addlegendentry{Normal, $\beta=\pgfmathprintnumber{\slopeStrong}$}

	\addplot[only marks, red, forget plot,error bars/.cd, y dir=both, y explicit] table [x={level}, y={variance}, y error={err}]{\WeakVar};
    \addplot[mark=none, red] table[x={level}, y={create col/linear regression={y=variance}}]{\WeakVar};
    \xdef\slopeWeak{\pgfplotstableregressiona}
	\addlegendentry{Binomial, $\beta=\pgfmathprintnumber{\slopeWeak}$}
		
\end{axis}

\begin{axis}[%
    name=plot2,
    at=(plot1.right of south east), anchor=left of south west,
    ymode = log,
    log basis y={2},
    width=0.45\textwidth,
    height=0.2\textheight,
    axis x line=bottom,
    axis y line=left,
    title=Convergence,
    xlabel=Level,
    ylabel={mean},
    legend pos= north east]

	\addplot[only marks, black, forget plot,error bars/.cd, y dir=both, y explicit] table[x={level}, y={mean}, y error={err}]{\BeloCon};
    \addplot[mark=none, black] table[x={level}, y={create col/linear regression={y=mean}}]{\BeloCon};
    \xdef\slopeBelomestny{\pgfplotstableregressiona}
	\addlegendentry{Euler, $p=\pgfmathprintnumber{\slopeBelomestny}$}

    \addplot[only marks, blue, forget plot,error bars/.cd, y dir=both, y explicit] table[x={level}, y={mean}, y error={err}]{\StrongCon};
    \addplot[mark=none, blue] table[x={level}, y={create col/linear regression={y=mean}}]{\StrongCon};
    \xdef\slopeStrong{\pgfplotstableregressiona}
	\addlegendentry{Normal, $p=\pgfmathprintnumber{\slopeStrong}$}

	\addplot[only marks, red,forget plot,error bars/.cd, y dir=both, y explicit] table [x={level}, y={mean}, y error=err]{\WeakCon};
    \addplot[mark=none, red] table[x={level}, y={create col/linear regression={y=mean}}]{\WeakCon};
    \xdef\slopeWeak{\pgfplotstableregressiona}
	\addlegendentry{Binomial, $p=\pgfmathprintnumber{\slopeWeak}$}

\end{axis}

\begin{axis}[%
    name=plot3,
    at=(plot1.below south east), anchor=above north east,
    ymode = log,
    xmode = log,
    log basis y={2},
    log basis x={2},
    width=0.45\textwidth,
    height=0.2\textheight,
    axis x line=bottom,
    axis y line=left,
    title={Paths on level $l=0$},
    xlabel={RMSE},
    ylabel={paths},
    legend pos = north east]

	\addplot[only marks, black, forget plot,error bars/.cd, y dir=both, y explicit] table[x={rmse}, y={paths}, y error={err}]{\BeloPath};
    \addplot[mark=none, black] table[x={rmse}, y={create col/linear regression={y=paths}}]{\BeloPath};
    \xdef\slopeBelomestny{\pgfplotstableregressiona}
	\addlegendentry{Euler, $slope=\pgfmathprintnumber{\slopeBelomestny}$}

    \addplot[only marks, blue, forget plot,error bars/.cd, y dir=both, y explicit] table[x={rmse}, y={paths}, y error={err}]{\StrongPath};
    \addplot[mark=none, blue] table[x={rmse}, y={create col/linear regression={y=paths}}]{\StrongPath};
    \xdef\slopeStrong{\pgfplotstableregressiona}
	\addlegendentry{Normal, $slope=\pgfmathprintnumber{\slopeStrong}$}

	\addplot[only marks, red,forget plot,error bars/.cd, y dir=both, y explicit] table [x={rmse}, y={paths}, y error={err}]{\WeakPath};
    \addplot[mark=none, red] table[x={rmse}, y={create col/linear regression={y=paths}}]{\WeakPath};
    \xdef\slopeWeak{\pgfplotstableregressiona}
	\addlegendentry{Binomial, $slope=\pgfmathprintnumber{\slopeWeak}$}
	
\end{axis}

\begin{axis}[%
    name=plot4,
    at=(plot3.right of south east), anchor=left of south west,
    width=0.45\textwidth,
    height=0.2\textheight,
    ymode = log,
    xmode = log,
    log basis y={2},
    log basis x={2},
    axis x line=bottom,
    axis y line=left,
    xlabel=$RMSE$,
    ylabel=$\mbox{RMSE}^2\mbox{cost})$,
    title=Computational Complexity,
    legend pos = south east]

\addplot[black, mark=*,error bars/.cd, y dir=both, y explicit] table [x={rmse}, y={cost}]{\BeloComp};
    \addlegendentry{Euler}

    \addplot[blue, mark=*,error bars/.cd, y dir=both, y explicit] table [x={rmse}, y={cost}]{\StrongComp};
    \addlegendentry{Normal}

    \addplot[red, mark=*,error bars/.cd, y dir=both, y explicit] table [x={rmse}, y={cost}, y error={err}]{\WeakComp};
    \addlegendentry{Binomial}

\end{axis}
\end{tikzpicture}
\caption{Computational results for the basket option.}\label{fig:EuropeanCall}
\end{figure}

We see that for the two Milstein schemes, the variance is  more or less identical,
and is approximately of order $\beta=2$, whereas the variance reduction for the Euler scheme is only of order $\beta =1$, as expected.
$\beta=2$ results by \cref{equ:compcost} in a quadratic computational complexity, which is also confirmed by the simulations.

Convergence is in all three cases as expected approximately of order 1, though for the Euler scheme slightly larger. Lastly, we see that the number of paths needed on level $l=0$ is significantly less when using enumeration than if we use standard Monte-Carlo. This then influences the computational complexity, such that the computational complexity of the Milstein scheme using enumeration is for $RMSE=1e-5$ about $4$ times smaller than the scheme without, though only about $2$ times for higher RMSEs due to the initial test simulation.

\section{Acknowledgements}
Nicky Cordua Mattsson would like to thank the SDU eScience center for partially funding his PhD.

Azadeh Ghasemifard would like to thank Dr.\ Mohammad Taghi Jahandideh, her supervisor in Isfahan University of Technology.
\def\bibsection{\section*{References}}
\bibliographystyle{siam}
\bibliography{bibliography}
\end{document}